\documentclass[12pt,reqno]{amsart}
\usepackage{enumerate, latexsym, amsmath, amsfonts, amssymb, amsthm, color}
\def\pmod #1{\ ({\rm{mod}}\ #1)}

\def\l{\left}
\def\r{\right}
\def\bg{\bigg}
\def\({\bg(}
\def\){\bg)}
\def\t{\text}
\def\f{\frac}

\def\bi{\binom}

\def\eq{\equiv}

\def\Proof{\noindent{\it Proof}}
\def\Ack{\medskip\noindent {\bf Acknowledgment}}
\theoremstyle{plain}
\newtheorem{theorem}{Theorem}

\newtheorem{lemma}{Lemma}

\theoremstyle{definition}

\theoremstyle{remark}
\newtheorem{remark}{Remark}

 \vspace{4mm}

\begin{document}

\hbox{Published on Ramanujan. J}
\medskip

\title
[{Proof of some congruence conjectures of Guo and Liu}]
{Proof of some congruence conjectures of Guo and Liu}

\author
[Guo-Shuai Mao] {Guo-Shuai Mao}

\address {(Guo-Shuai Mao) Department of Mathematics, Nanjing
University, Nanjing 210093, People's Republic of China}
\email{mg1421007@smail.nju.edu.cn}
\keywords{Central binomial coefficients; Congruences; Bernoulli numbers; Zeilberger algorithm.
\newline \indent This research was supported by the National Natural Science Foundation (Grant No. 11571162) of China.}
\begin{abstract} Let $n$ and $r$ be positive integers. Define the numbers $S_n^{(r)}$ by $S_n^{(r)}=\sum_{k=0}^n\binom{n}{k}^2\binom{2k}{k}(2k+1)^r.$
In this paper we prove some conjectures of Guo and Liu which extend some conjectures of Z.-W. Sun \cite{Su1}, such as:
There exist integers $a_{2r-1}$ and $b_r$, independent of $n$, such that
$$a_{2r-1}\sum_{k=0}^{n-1}S_k^{(2r-1)}\equiv0\pmod{n^2}\ \mbox{and}\ b_r\sum_{k=0}^{n-1}kS_k^{(r)}\equiv0\pmod{n^2}.$$
By Zeilberger algorithm, we find that for all $0\leq j<n$, $$(2j+1)\binom{2j}j\sum_{k=j}^{n-1}(2k-j+1)\binom kj^2\equiv0\pmod{n^2}.$$

\end{abstract}
\maketitle
\section{Introduction}
\setcounter{lemma}{0}
\setcounter{theorem}{0}
\setcounter{corollary}{0}
\setcounter{remark}{0}
\setcounter{equation}{0}
\setcounter{conjecture}{0}
Differential equations of the Calabi-Yau type \cite{AESZ} have the form
$$D_{n,k}y:=\left(\sum_{i=0}^kP_i(\theta)\right)y=0,$$
where $\theta:=z\frac{\partial}{\partial z}$ and the $P_i(\theta)$ are polynomials of degree $n$ with integer coefficients.
Such equations are related to Ap\'{e}ry's proof of the irrationality of $\zeta(3)$. Ap\'{e}ry's proof involved the sequence $b_n=\sum_{k=0}^n\binom{n}{k}^2\binom{n+k}{k}^2$. Zagier \cite{Z} has made a study of similar sequences. One such sequence is $g_n=\sum_{k=0}^n\binom{n}{k}^2\binom{2k}k$. The function $y_0(z)=\sum_{n=0}^{\infty}A_nz^n$, where $A_n=\binom{2n}{n}^2g_n$ is the analytic solution of the $4$-th order differential equation
$$Dy=(\theta^4-4z(2\theta+1)^2(10\theta^2+10\theta+3)+144z^2(2\theta+1)^2(2\theta+3)^2)y=0,$$
which satisfies the condition $y_0(0)=A_0=1$.

The author and Z.-W. Sun \cite{MS} made a study of the sequence $g_n$.

Z.-W. Sun \cite{Su1} defined
$$S_n=\sum_{k=0}^n\binom{n}k^2\binom{2k}{k}(2k+1)\ (n=0,1,2,\ldots),$$
and gave a nice property of the sum $\sum_{k=0}^{n-1}S_k$ \cite[(1.19)]{Su1}:
$$\frac1{n^2}\sum_{k=0}^{n-1}S_k=\sum_{k=0}^{n-1}\binom{n-1}k^2C_k\in \mathbb{Z},$$
where $C_k=\frac1{k+1}\binom{2k}k$ is the $k$-th Catalan number.\\

 Guo and Liu \cite{GL} confirmed some conjectures of Z.-W. Sun for $S_n$ and gave an identity of the sum involving $kS_k$ \cite[(2.1)]{GL}
$$\frac{4}{n^2}\sum_{k=0}^{n-1}kS_k=\sum_{k=0}^{n-1}\l(6k\binom{n-1}k^2+\binom{n-1}k\binom{n-1}{k+1}\r)C_k.$$
 Z.-W. Sun also defined $S_n^{-}=\sum_{k=0}^n\binom{n}{k}^2\binom{2k}{k}(2k+1)^2(-1)^k$, and gave some conjectures for these numbers \cite[Conjecture 5.6]{Su1}. Guo and Liu \cite{GL} proved one of those conjectures and defined
$$S_n^{(r)}=\sum_{k=0}^n\binom{n}{k}^2\binom{2k}{k}(2k+1)^r\  \mbox{,}\  T_n^{(r)}=\sum_{k=0}^n\binom{n}{k}^2\binom{2k}{k}(2k+1)^r(-1)^k,$$ which extend $S_n$ and $S_n^{-}$.

Guo and Liu \cite{GL1} confirmed a conjecture of Z.-W. Sun \cite{Su1} involving $R_n=\sum_{k=0}^n\binom{n+k}{2k}\binom{2k}{k}\frac{1}{2k-1}$.

The following theorem confirms the Conjecture 5.1 of \cite{GL}.
\begin{theorem}\label{Th1.1} Let $n$ and $r$ be positive integers and $p$  be a prime, then we have
\begin{equation}\label{1.2}
\sum_{k=0}^{n-1}S_k^{(2r)}\equiv0\pmod{n^2},
\end{equation}
\begin{equation}\label{1.3}
\sum_{k=0}^{n-1}T_k^{(2r)}\equiv0\pmod{n^2},
\end{equation}
and
\begin{equation}\label{1.4}
\sum_{k=0}^{p-1}T_k^{(2)}\equiv\frac{p^2}2\l(5-3\(\frac{5}p\)\r)\pmod{p^3}.
\end{equation}
\end{theorem}
In \cite{KM} we know $S_m(n)=1^m+2^m+\cdots+(n-1)^m$ and the m-th Bernoulli number can be written as $B_m=U_m/V_m$  where $(U_m,V_m)=1$ and $V_m>0$.
The Bernoulli numbers are given by
$$\f x{e^x-1}=\sum_{n=0}^\infty B_n\f{x^n}{n!}\ \ (0<|x|<2\pi).$$
We confirm the Conjecture 5.2 of Guo and Liu \cite{GL} by the following theorem.
\begin{theorem}\label{Th1.2} Let $n$ and $r$ be positive integers. Let
\begin{equation*}
a_{2r-1}=\begin{cases}1&\t{if}\ \t{r=1},\\ \frac12V_{2r-2}&\t{if}\ \t{r}{>}1.\end{cases}
\end{equation*}
and
\begin{equation*}
b_r=\begin{cases}2V_{r-1}&\t{if}\ \t{r}\eq1\pmod2,\\2V_r&\t{if}\ \t{r}\eq0\pmod2.\end{cases}
 \end{equation*}
Then
\begin{equation}\label{1.5}
a_{2r-1}\sum_{k=0}^{n-1}S_k^{(2r-1)}\equiv0\pmod{n^2},
\end{equation}
\begin{equation}\label{1.6}
b_r\sum_{k=0}^{n-1}kS_k^{(r)}\equiv0\pmod{n^2}.
\end{equation}
\end{theorem}
\begin{remark} By the formulas of $a_{2r-1}$ and $b_r$ we can partially answer conjecture 5.3 of Guo and Liu \cite{GL}.

First we give some Lemmas in Section 2. Then we prove Theorems \ref{Th1.1} and \ref{Th1.2} in Section 3 and Section 4 respectively.
\end{remark}
\section{Some Preliminaries}
\setcounter{lemma}{0}
\setcounter{theorem}{0}
\setcounter{corollary}{0}
\setcounter{remark}{0}
\setcounter{equation}{0}
\setcounter{conjecture}{0}
In order to prove the above theorems, we give some lemmas first:

\begin{lemma}\label{Lem1.1}\cite[(1.48) and (6.30)]{G} Let $m$ and $n$ be nonnegative integers, then
 \begin{equation}\label{1.7}
 \sum_{k=0}^n\binom{x+k}{m}=\binom{n+x+1}{m+1}-\binom{x}{m+1},
 \end{equation}
 and
 \begin{equation}\label{1.8}\sum_{k=0}^n\binom{n}{k}^2\binom{x+k}{2n}=\binom{x}{n}^2.
 \end{equation}

 \end{lemma}

 \begin{lemma}\label{Lem1.2}{\cite[Lemma 3.2]{MS}} For any nonnegative integer $n$, we have
\begin{equation}\label{1.9}
 \sum_{k=0}^n\binom{n}{k}^2\binom{x+k}{2n+1}=\f1{(4n+2)\bi{2n}n}\sum_{k=0}^n(2x-3k)\bi xk^2\bi{2k}k.
 \end{equation}
 \end{lemma}
The identity (\ref{1.9}) is crucial to our proofs of Theorems \ref{Th1.1}-\ref{Th1.2}.

We know that the Euler numbers $\{E_n\}$ and Euler polynomials $\{E_n(x)\}$ are defined by
$$\frac{2e^t}{e^{2t}+1}=\sum_{n=0}^\infty E_n\frac{t^n}{n!}\ (|t|<\frac{\pi}2)\ \mbox{and}
\ \frac{2e^{xt}}{e^{t}+1}=\sum_{n=0}^\infty E_n(x)\frac{t^n}{n!}\ (|t|<\pi).$$
It is known that the Euler numbers and Euler polynomials have some properties, such as $E_{2n-1}=0, E_n\in \mathbb{Z}$ for all positive integer $n$, $E_n(x)+E_n(x+1)=2x^n$, $E_n(\frac12)=E_n/2^n$ and
\begin{align}\label{1.1}
E_n(x+y)=\sum_{k=0}^n\binom{n}kE_k(x)y^{n-k}.
\end{align}
\begin{lemma}\label{Lem2.1} Let $n$ and $r$ be positive integers, then we have
$$\sum_{k=0}^{n-1}(2k+1)^{2r-1}\equiv 0\pmod n.$$
\end{lemma}
\Proof:  We know that $$\sum_{k=0}^{n-1}(2k+1)^{2r-1}=\sum_{j=0}^{n-1}(2n-2j-1)^{2r-1}\equiv-\sum_{j=0}^{n-1}(2j+1)^{2r-1}\pmod{2n},$$
hence
$$2\sum_{k=0}^{n-1}(2k+1)^{2r-1}\equiv0\pmod{2n}.$$
Therefore
 $$\sum_{k=0}^{n-1}(2k+1)^{2r-1}\equiv0\pmod n.$$
This confirms Lemma \ref{Lem2.1}. \qed

\begin{lemma}\label{Lem2.2} For any positive integer $n$ and $r$, we have
$$\sum_{k=0}^{n-1}(2k+1)^{2r-1}(-1)^k\equiv 0\pmod n.$$
\end{lemma}
\Proof: Combining $E_n(x)+E_n(x+1)=2x^n$ with $E_n(\frac12)=E_n/2^n$ and (\ref{1.1}), we have
\begin{align*}
 &\sum_{k=0}^{n-1}(-1)^k(2k+1)^{2r-1}
 \\=&2^{2r-2}\sum_{k=0}^{n-1}\left((-1)^kE_{2r-1}\left(k+\frac12\right)-(-1)^{k+1}E_{2r-1}\left(k+\frac32\right)\right)
 \\=&2^{2r-2}\left(E_{2r-1}\left(\frac12\right)-(-1)^nE_{2r-1}\left(n+\frac12\right)\right)
 \\=&\frac{E_{2r-1}}2-(-1)^n\frac{E_{2r-1}}2-(-1)^n2^{2r-2}\sum_{k=0}^{2r-2}\binom{2r-1}k\frac{E_k}{2^k}n^{2r-1-k}.
\end{align*}
Note that $E_{2n-1}=0$, $E_n\in\mathbb{Z}$ for all $n\geq1$ and $2^{2r-2}/2^k\in\mathbb{Z}$ for all $k=0,1,\cdots,2r-2$. So we obtain the desired result. \qed
\section{Proof of Theorem 1.1}
\setcounter{lemma}{0}
\setcounter{theorem}{0}
\setcounter{corollary}{0}
\setcounter{remark}{0}
\setcounter{equation}{0}
\setcounter{conjecture}{0}
{\it Proof of \rm(\ref{1.2})}: By (\ref{1.8}) and (\ref{1.7}) we obtain that
$$\sum_{k=0}^{n-1}S_k^{(2r)}=\sum_{j=0}^{n-1}\binom{2j}j(2j+1)^{2r}\sum_{l=0}^j\binom{j}l^2\binom{n+l}{2j+1}.$$
By using Lemma \ref{Lem1.2} we deduce that
$$\sum_{k=0}^{n-1}S_k^{(2r)}=\frac12\sum_{j=0}^{n-1}(2j+1)^{2r-1}\sum_{l=0}^j(2n-3l)\binom{n}l^2\binom{2l}l.$$
With the help of Lemma \ref{Lem2.1} we have
\begin{align*}
&\sum_{k=0}^{n-1}S_k^{(2r)}
\\\equiv&\sum_{l=1}^{n-1}\frac{n^2(2n-3l)}{l^2}\binom{n-1}{l-1}^2\binom{2l-1}{l-1}\sum_{j=l}^{n-1}(2j+1)^{2r-1}\pmod{n^2}.
\end{align*}
Note that
\begin{align}\label{2.0}
\frac{2n-3l}l\binom{n-1}{l-1}=2\binom{n}l-3\binom{n-1}{l-1}\in\mathbb{Z},\ \frac nl\binom{n-1}{l-1}=\binom nl\in\mathbb{Z}
\end{align}
and by using Lemma \ref{Lem2.1} again we have
\begin{align*}
\frac{n^2(2n-3l)}{l^2}\binom{n-1}{l-1}^2\sum_{j=0}^{n-1}(2j+1)^{2r-1}\equiv0\pmod{n^2}
\end{align*}
and
\begin{align*}
\frac{n^2(2n-3l)}{l^2}\binom{n-1}{l-1}^2\sum_{j=0}^{l-1}(2j+1)^{2r-1}\equiv0\pmod{n^2}.
\end{align*}
Hence
\begin{equation*}
\frac{n^2(2n-3l)}{l^2}\binom{n-1}{l-1}^2\sum_{j=l}^{n-1}(2j+1)^{2r-1}\eq0\pmod{n^2}.
\end{equation*}
Therefore  $$\sum_{k=0}^{n-1}S_k^{(2r)}\equiv0\pmod{n^2}.$$
This confirms (\ref{1.2}). \qed

\noindent{\it Proof of \rm(\ref{1.3})}:\ By (\ref{1.8}) and (\ref{1.7}) we have
$$\sum_{k=0}^{n-1}T_k^{(2r)}=\sum_{j=0}^{n-1}\binom{2j}j(2j+1)^{2r}(-1)^j\sum_{l=0}^j\binom{j}l^2\binom{n+l}{2j+1}.$$
This, together with Lemma \ref{Lem1.2} yield that
\begin{equation}\label{2.1}
\sum_{k=0}^{n-1}T_k^{(2r)}=\frac12\sum_{j=0}^{n-1}(2j+1)^{2r-1}(-1)^j\sum_{l=0}^j(2n-3l)\binom{n}l^2\binom{2l}l.
\end{equation}
By Lemma \ref{Lem2.2} we have
\begin{align*}
&\sum_{k=0}^{n-1}T_k^{(2r)}
\\\equiv&\sum_{l=1}^{n-1}\frac{n^2(2n-3l)}{l^2}\binom{n-1}{l-1}^2\binom{2l-1}{l-1}\sum_{j=l}^{n-1}(2j+1)^{2r-1}(-1)^j\pmod{n^2}.
\end{align*}
By (\ref{2.0}) and Lemma \ref{Lem2.2} again we have
\begin{align*}
\frac{n^2(2n-3l)}{l^2}\binom{n-1}{l-1}^2\sum_{j=0}^{n-1}(2j+1)^{2r-1}(-1)^j\equiv0\pmod{n^2}
\end{align*}
and
\begin{align*}
\frac{n^2(2n-3l)}{l^2}\binom{n-1}{l-1}^2\sum_{j=0}^{l-1}(2j+1)^{2r-1}(-1)^j\equiv0\pmod{n^2}.
\end{align*}
Therefore
$$\frac{n^2(2n-3l)}{l^2}\binom{n-1}{l-1}^2\sum_{j=l}^{n-1}(2j+1)^{2r-1}(-1)^j\equiv0\pmod{n^2}.$$
Hence
$$\sum_{k=0}^{n-1}T_k^{(2r)}\equiv0\pmod{n^2}.$$
Thus we have completed the proof of (\ref{1.3}). \qed

\noindent{\it Proof of \rm(\ref{1.4})}: Let $n=p$ and $r=1$ in (\ref{2.1}) we have
$$\sum_{k=0}^{p-1}T_k^{(2)}=\frac12\sum_{j=0}^{p-1}(2j+1)(-1)^j\sum_{l=0}^j(2p-3l)\binom{p}l^2\binom{2l}l.$$
By $\sum_{j=0}^{p-1}(2j+1)(-1)^j=p$, $\sum_{j=l}^{p-1}(2j+1)(-1)^j=p+(-1)^ll$ and $$\binom{p-1}{l-1}^2=\prod_{j=1}^{l-1}\l(\frac pj-1\r)^2\equiv1\pmod p$$ we obtain that
$$\sum_{k=0}^{p-1}T_k^{(2)}\equiv p^2-\frac32p^2\sum_{l=1}^{p-1}\binom{2l}l(-1)^l\pmod{p^3}.$$
It is known that $\binom{2l}l=\binom{-\frac12}l(-4)^l\equiv\binom{\frac{p-1}2}l(-4)^l\pmod p$, therefore
\begin{align*}
&\sum_{k=0}^{p-1}T_k^{(2)}\equiv p^2-\frac32p^2\sum_{l=1}^{(p-1)/2}\binom{\frac{p-1}2}l4^l
\\ \equiv&p^2-\frac32p^2\l(5^{(p-1)/2}-1\r)\equiv\frac{p^2}2\l(5-3\(\frac5p\)\r)\pmod {p^3}.
\end{align*}
Obviously $\big(\frac5p\big)=\big(\frac p5\big)$, so (\ref{1.4}) follows. \qed
\section{Proof of Theorem 1.2}
\setcounter{lemma}{0}
\setcounter{theorem}{0}
\setcounter{corollary}{0}
\setcounter{remark}{0}
\setcounter{equation}{0}
\setcounter{conjecture}{0}
Recall that the m-th Bernoulli number can be written as $B_m=U_m/V_m$  where $(U_m,V_m)=1$ and $V_m>0$.
\begin{lemma}\label{Lem3.1}\cite[Proposition 15.2.2]{KM}  If $m\geq2$ is an even integer, then for all $n\geq1$ we have
$$V_mS_m(n)\equiv nU_m\pmod{n^2}.$$
\end{lemma}
\begin{lemma}\label{Lem3.2} Let $n$ and $r$ be positive integers and $r>1$, then
$$V_{2r-2}\sum_{j=0}^{n-1}(2j+1)^{2r-2}\equiv0\pmod {2n}.$$
\end{lemma}
\Proof: We know $\sum_{j=0}^{n-1}(2j+1)^{2r-2}=S_{2r-2}(2n)-2^{2r-2}S_{2r-2}(n)$, so that by Lemma \ref{Lem3.1} we obtain that
$$V_{2r-2}S_{2r-2}(2n)\equiv2nU_{2r-2}\equiv0\pmod {2n}$$ and $$V_{2r-2}S_{2r-2}(n)\equiv nU_{2r-2}\equiv0\pmod n$$ for all $r>1$. So for all $r>1$ we have $$V_{2r-2}\sum_{j=0}^{n-1}(2j+1)^{2r-2}\equiv0\pmod {2n}.$$
This confirms Lemma \ref{Lem3.2}. \qed

\noindent{\it Proof of \rm(\ref{1.5})}: For $r=1$ we know $\sum_{k=0}^{n-1}S_k\equiv0\pmod{n^2}$, which was given by Z.-W. Sun, so we just need to show the result for all $r>1$. By (\ref{1.8}) and (\ref{1.7}) we have
$$\sum_{k=0}^{n-1}S_k^{(2r-1)}=\sum_{j=0}^{n-1}\binom{2j}j(2j+1)^{2r-1}\sum_{l=0}^j\binom{j}l^2\binom{n+l}{2j+1}.$$
Then via Lemma \ref{Lem1.2} we obtain
$$\sum_{k=0}^{n-1}S_k^{(2r-1)}=\frac{1}2\sum_{j=0}^{n-1}(2j+1)^{2r-2}\sum_{l=0}^j(2n-3l)\binom{n}l^2\binom{2l}l.$$
Therefore
\begin{align*}
&V_{2r-2}\sum_{k=0}^{n-1}S_k^{(2r-1)}-V_{2r-2}n\sum_{j=0}^{n-1}(2j+1)^{2r-2}
\\=&\frac{V_{2r-2}}2\sum_{j=1}^{n-1}(2j+1)^{2r-2}\sum_{l=1}^j(2n-3l)\binom{n}l^2\binom{2l}l
\\=&\sum_{l=1}^{n-1}\frac{V_{2r-2}n^2(2n-3l)}{l^2}\binom{n-1}{l-1}^2\binom{2l-1}{l-1}\sum_{j=l}^{n-1}(2j+1)^{2r-2}.
\end{align*}
This, together with Lemma \ref{Lem3.2} yield that
\begin{align*}
&V_{2r-2}\sum_{k=0}^{n-1}S_k^{(2r-1)}
\\\equiv &\sum_{l=1}^{n-1}\frac{V_{2r-2}n^2(2n-3l)}{l^2}\binom{n-1}{l-1}^2\binom{2l-1}{l-1}\sum_{j=l}^{n-1}(2j+1)^{2r-2}\pmod{2n^2}.
\end{align*}
By (\ref{2.0}) and Lemma \ref{Lem3.2} we have
\begin{align*}
\frac{V_{2r-2}n^2(2n-3l)}{l^2}\binom{n-1}{l-1}^2\sum_{j=0}^{n-1}(2j+1)^{2r-2}\equiv0\pmod{2n^2}
\end{align*}
and
\begin{align*}
\frac{V_{2r-2}n^2(2n-3l)}{l^2}\binom{n-1}{l-1}^2\sum_{j=0}^{l-1}(2j+1)^{2r-2}\equiv0\pmod{2n^2}.
\end{align*}
Hence
$$\frac{V_{2r-2}n^2(2n-3l)}{l^2}\binom{n-1}{l-1}^2\sum_{j=l}^{n-1}(2j+1)^{2r-2}\equiv0\pmod{2n^2}.$$
Therefore $$V_{2r-2}\sum_{k=0}^{n-1}S_k^{(2r-1)}\equiv0\pmod {2n^2},$$
so we just need to take
\begin{equation*}
a_{2r-1}=\begin{cases}1&\t{if}\ \t{r=1},\\ \frac12V_{2r-2}&\t{if}\ \t{r}{>}1.\end{cases}
\end{equation*}
So (\ref{1.5}) follows. \qed
\begin{lemma}\label{Lem3.3} For any positive integer $n$ and $0\leq j<n$ we have
$$(2j+1)\binom{2j}j\sum_{k=j}^{n-1}(2k-j+1)\binom kj^2\equiv0\pmod{n^2}.$$
\end{lemma}
\Proof: Set $u_j=(2j+1)\binom{2j}j\sum_{k=j}^{n-1}(2k-j+1)\binom kj^2$. By Zeilberger \cite{PWZ} we obtain the explicit formula
$$u_j=\frac{(1+2j)n^2}{j+1}\binom{2j}j\binom{n-1}j^2$$ for all $0\leq j<n-1$.
It is known that $\frac{1}{j+1}\binom{2j}j=C_j\in\mathbb{Z}$, so we have $u_j\equiv0\pmod {n^2}$ for all $0\leq j<n-1$ and $$u_{n-1}=(2n-1)\binom{2n-2}{n-1}n=n^2\binom{2n-1}{n-1}\equiv0\pmod {n^2}.$$ Hence $u_j\equiv0\pmod {n^2}$ for all $0\leq j<n$.\qed

\noindent{\it Proof of \rm(\ref{1.6})}: First we have
\begin{align*}
&\sum_{k=0}^{n-1}\left(4kS_k^{(r)}-S_k^{(r+1)}+3S_k^{(r)}\right)
\\=&\sum_{k=0}^{n-1}\sum_{j=0}^k(4k-2j+2)\binom kj^2\binom{2j}j(2j+1)^r
\\=&\sum_{j=0}^{n-1}\binom{2j}j(2j+1)^r\sum_{k=j}^{n-1}(4k-2j+2)\binom kj^2=2\sum_{j=0}^{n-1}(2j+1)^{r-1}u_j.
\end{align*}
By Lemma \ref{Lem3.3} we deduce that $$\sum_{k=0}^{n-1}\left(4kS_k^{(r)}-S_k^{(r+1)}+3S_k^{(r)}\right)\equiv0\pmod{2n^2}.$$
We know when $r=1$, $\sum_{k=0}^{n-1}4kS_k^{(r)}\equiv0\pmod{n^2}$, which was given by Guo and Liu \cite{GL}. So we just need to show that for all $r>1$ .

\noindent{\it \rm{Case 1}}. If $2|r$, then by (\ref{1.2}) we have $V_r\sum_{k=0}^{n-1}S_k^{(r)}\equiv0\pmod {2n^2}$ since 6 always divides the denominator of $B_r$ \cite[(p.233)]{KM}. By (\ref{1.5}) we have $V_r\sum_{k=0}^{n-1}S_k^{(r+1)}\equiv0\pmod{2n^2}$. Hence $$V_r\sum_{k=0}^{n-1}4kS_k^{(r)}\equiv0\pmod{2n^2}.$$
{\it \rm{Case 2}}. If $2\nmid r$ and $r>1$, then by (\ref{1.2}) we have $V_{r-1}\sum_{k=0}^{n-1}S_k^{(r+1)}\equiv0\pmod {2n^2}$ since 6 always divides the denominator of $B_{r-1}$. By (\ref{1.5}) we have $V_{r-1}\sum_{k=0}^{n-1}S_k^{(r)}\equiv0\pmod{2n^2}$. Thus we have $$V_{r-1}\sum_{k=0}^{n-1}4kS_k^{(r)}\equiv0\pmod{2n^2}.$$
So combining  Case 1 and Case 2 we finally obtain that there exist integers $b_r$ such that $$b_r\sum_{k=0}^{n-1}kS_k^{(r)}\equiv0\pmod{n^2}.$$
Hence we just need to take
\begin{equation*}
 b_r=\begin{cases}2V_{r-1}&\t{if}\ \t{r}\eq1\pmod2,\\2V_r&\t{if}\ \t{r}\eq0\pmod2.\end{cases}
\end{equation*}
So (\ref{1.6}) follows. \qed

\Ack. The author would like to thank Prof. Z.-W. Sun, Prof. Hao Pan and the referees for helpful comments.


\begin{thebibliography}{AESZ}
\bibitem[AESZ]{AESZ} G. Almkvist, C. van. Enckevort, D. van. Straten, W. Zudilin, {\it Tables of Calabi-Yau equations}, Preprint, {\tt arXiv:math/0507430}, 9 October 2010.
 \bibitem[G]{G}     H.W. Gould, {\it Combinatorial Identities}, Morgantown Printing and Binding Co., 1972.
\bibitem[GL]{GL}    V.J.W. Guo and J.-C. Liu, {\it Proof of some conjectures of Z.-W.Sun on the divisibility of certain double-sums}, Int. J. Number Theory {\bf 12} (2016), 615--623.
\bibitem[GL1]{GL1}  V.J.W. Guo and J.-C. Liu, {\it Proof of a conjecture of Z.-W.Sun on the divisibility of a triple-sum}, {\bf 156} (2015), 154--160.
\bibitem[KM]{KM}    K. Ireland and M. Rosen, {\it A Classical Introduction to Modern Number Theory}, 2nd Editor, Grad. Texts in Math., Vol.84, Springer, New York, 1990.
\bibitem[MS]{MS}    G.-S. Mao and Z.-W. Sun, {\it Two congruences involving harmonic numbers with applications}, Int. J. Number Theory {\bf 12} (2016), no.02, 527-539.
\bibitem[PWZ]{PWZ}  M.Petkov\u{s}ek, H.S. Wilf and D. Zeilberger, {\it $A=B$}, A K Peters, Ltd., Wellesley, MA, 1996.
\bibitem[Su1]{Su1}  Z.-W. Sun, {\it Two new kinds of numbers and related divisibility results}, Colloq. Math. {\bf 154} (2018), no.2, 241--273..
\bibitem[Z]{Z}      D. Zagier, {\it Integral solutions of Ap\'{e}ry-like recurrence equations}, Groups and symmetries,
CRM Proc. Lecture Notes, vol. 47, Amer. Math. Soc., Providence, RI, 2009, pp. 349--366. MR 2500571
\end{thebibliography}
\end{document}